\documentclass[12pt,leqno]{article}

\newcommand\qed{\hfill$\sqcap\kern-8.0pt\hbox{$\sqcup$}$}
\newcommand\NN{\hbox{I\kern-.2em\hbox{N}}}
\newcommand\RR{\hbox{I\kern-.2em\hbox{R}}}
\newcommand\sRR{{\sl \hbox{I\kern-.2em\hbox{R}}}}
\newcommand{\PP}{{\bf P}^k}

\usepackage{amsmath}
\usepackage{amssymb}
\usepackage{latexsym}
\usepackage{amscd}
\usepackage[latin1]{inputenc}
\usepackage[french]{babel}

\newcommand{\fin}{\hfill{\Large$\Box$}}
\newtheorem{theo}{Th\'eor\`eme}
\newtheorem{prop}{Proposition}
\newtheorem{lem}{Lemme}
\newtheorem{rem}{Remarque}
\newtheorem{cor}{Corollaire}
\newtheorem{defi}{D\'efinition}

\newcommand\ZZ{{{\rm Z}\kern-.28em{\rm Z}}}
\newcommand\QQ{\hbox{I\kern-.53em\hbox{Q}}}
\newcommand{\Id}{{\rm Id_{{\bf C}^k}}}

\newcommand{\Lip}{{\rm Lip \, }}
\newcommand{\Jac}{{\rm Jac \,  }}

\begin{document}

\date{}

\title 
{Une caractérisation des endomorphismes de Lattès par leur mesure
de Green}

\vskip0.5cm
\author { F. Berteloot et C. Dupont}
\vskip0.5cm
\maketitle
\noindent{ABSTRACT.}\\
{\small We show that the Latt\` es endomorphisms are the only holomorphic endomorphisms
of the complex $k$-dimensional projective space whose measure of maximal
entropy is absolutely continuous with respect to the Lebesgue measure. As a
consequence, Latt\` es endomorphisms are also characterized by other extremal
properties as the maximality of the Hausdorff dimension of their measure of
maximal entropy or the minimality of their Liapounov exponents. Our proof uses
a linearization method which is of independant interest and a previous characterization by
the regularity of the Green current.
}
\noindent{AMS-CLASSIFICATION:} {32H50, 32U40, 37C45.
}\\

\section{Introduction et résultats}
Les propriétés dynamiques d'un endomorphisme holomorphe $f$ de degré algébrique $d\ge 2$ sur l'espace projectif
complexe ${\bf P}^k$ se reflètent sur son courant et sa mesure de Green. Le courant de Green, noté $T$, est un $(1,1)$-courant positif fermé, obtenu comme la limite des $(1,1)$-formes $\frac{1}{d^n}f^{n*} \omega$, où $\omega$
désigne la forme de Fubini-Study. La mesure de Green, notée $\mu$, est une mesure de probabilité invariante, obtenue comme $k$-ième
puissance extérieure de $T$. Ces objets, introduits par Hubbard-Papadopol \cite{HP} et Fornaess-Sibony
\cite{FS}, possèdent de remarquables propriétés ergodiques. Fornaess et Sibony ont montré
que la mesure de Green est mélangeante \cite{FS2}. Briend et Duval ont établi que ses exposants de
Liapounov sont supérieurs à $\log \sqrt{d}$  \cite{BD} et qu'elle est l'unique mesure d'entropie maximale de $f$ \cite{BD2}.\\

La dimension de $\mu$, notée $\dim (\mu)$, est définie comme la borne inférieure des dimensions de
Hausdorff des boréliens de $\mu$-mesure pleine. C'est une caractéristique géométrique importante du système
dynamique $({\bf P}^k, f, \mu)$. L'une des premières questions concernant l'estimation 
de cette dimension est de déterminer les
systèmes pour lesquels elle est maximale ou, ce qui s'avère équivalent, ceux
dont la mesure de Green est absolument continue 
par rapport à la mesure de Lebesgue $\omega^k$. L'objet de cet article est de caractériser ces systèmes,
ce qui répond à une question posée par Fornaess et Sibony dans \cite{FS4} :

\begin{theo}\label{TH} 
Les seuls endomorphismes holomorphes de ${\bf P}^k$ dont la mesure de Green est absolument continue par
rapport à $\omega^k$ sont les endomorphismes de Lattès.
\end{theo}
 
Rappelons qu'un endomorphisme $f$ est de Lattès si il fait commuter un diagramme : 
\begin{equation*}
\begin{CD}
{\bf C}^k @ > D >>                {\bf C}^k\\
{\sigma}  @ VVV                     @VVV{\sigma}\\
{\bf P}^k @ > f >>                {\bf P}^k
\end{CD}
\end{equation*}    
où $D$ est une application affine de partie
lin\'eaire $\sqrt d\,U$ (où $U$ est unitaire) et $\sigma$ un rev\^etement ramifi\'e
sur les fibres duquel un groupe cristallographique complexe  
agit transitivement. De tels endomorphismes existent en toute dimension $k$ et tout degr\'e $d$; 
leur mesure
de Green est absolument continue par rapport à $\omega^k$ d'exposants de Liapounov 
égaux à $\log \sqrt d$ \cite{Du1}. En dimension $1$, ils co\"incident avec les 
fractions rationnelles induites par une isogénie d'un tore complexe au moyen d'une fonction 
elliptique. Ils sont traditionnellement appelés "exemples de
Lattès" et font l'objet d'une étude détaillée dans l'article de revue de Milnor \cite{Mi}.  
Signalons enfin que les endomorphismes de Lattès interviennent naturellement dans d'autres 
problèmes, comme celui de la densité des fractions rationnelles hyperboliques via la "conjecture
NILF"  
(voir \cite{MM}, \cite{BM} Chap. 7) ou celui de la classification des paires d'endomorphismes qui
 commutent \cite{DS}. 
Ils fournissent également des exemples surprenants de domaines de ${\bf C}^{k+1}$ 
munis d'auto-applications holomorphes propres non injectives \cite{Du1}.\\

Voyons comment le théorème \ref{TH} se traduit en terme de dimension de la mesure. 
La maximalité de la dimension entra\^\i ne la minimalité des exposants. 
Cela résulte de l'inégalité $\dim(\mu) \leq 2(k-1) + {\log d \over \lambda_k}$, où $\lambda_k$ désigne 
le plus grand exposant de $\mu$. Cette estimation, dont la preuve est esquissée en appendice, 
est due à Binder et 
DeMarco \cite{BdM} pour les applications polynomiales (voir aussi \cite{DiDu} pour un résultat 
plus précis). Il est alors possible d'adapter aux dimensions 
supérieures le travail de Ledrappier  \cite{L1}, \cite{L2} selon lequel, pour toute fraction 
rationnelle, l'égalité dans la formule de Margulis-Ruelle 
entra\^\i ne l'absolue continuité de $\mu$. Cela fait l'objet de \cite{D} et concerne en particulier les 
mesures de Green d'exposants minimaux. Le théorème \ref{TH} admet donc pour corollaire :

\begin{cor}\label{COR} 
Soit un système $({\bf P}^k, f, \mu)$ où $f$ est de degr\'e $d$. Les propri\'et\'es suivantes sont
équivalentes :
\begin{itemize}
\item[1.] La dimension de $\mu$ est maximale, \'egale à $2k$.
\item[2.] Les exposants de $\mu$ sont minimaux, égaux à $\log \sqrt{d}$.
\item[3.] L'endomorphisme $f$ est de Lattès.
\end{itemize}
\end {cor}

Ainsi, pour un système $({\bf P}^k, f, \mu)$ générique,  la mesure $\mu$ est singulière par rapport à $\omega^k$, l'un de ses exposants est strictement supérieur à $\log \sqrt{d}$ et sa dimension est strictement inférieure à $2k$. \\   

En dimension $k=1$, on trouve une démonstration du théorème \ref{TH} dans l'article de Mayer \cite{M}. 
Signalons aussi le résultat beaucoup plus précis de Zdunik \cite {Z} qui stipule que la dimension de 
$\mu$ coïncide avec celle de son support (l'ensemble de Julia de $f$) si et seulement
si $f$ est un exemple de Lattès, un polyn\^ ome de Tchebychev ou une puissance $z^{\pm d}$.
La démonstration de Mayer repose sur un procédé de linéarisation  
consistant à comparer les itérées $f^n$ avec leurs applications linéaires tangentes $d_x f^n$. 
Un tel procédé permet de ``régulariser'' la densité mesurable de $\mu$ : celle-ci est en fait lisse sur 
un ouvert. La structure de $f$ se lit alors sur l'équation fonctionnelle $f^* \mu =d \mu$.  \\                   

Il y a plusieurs difficultés à surmonter en dimension supérieure.
Fondamentalement, le probl\` eme  tient \` a ce que la mesure $\mu$ ne porte pas les informations 
géométriques 
"directionnelles" nécessaires à l'analyse de la structure de $f$ : celles-ci sont recel\'ees par le 
courant $T$ dont elle dérive ($\mu=T^k$) et s'y lisent particulièrement bien lorsque celui-ci est lisse :\\ 

{\bf Théorème} (Berteloot-Loeb \cite{BL})
{\it Tout endomorphisme holomorphe de ${\bf P}^k$ dont le courant de Green coïncide 
avec une $(1,1)$-forme lisse strictement positive sur un ouvert est un exemple de Lattès.}\\ 

\noindent Il s'agit donc de déduire la régularité du courant $T$ de l'absolue continuité de la
mesure $\mu = T^k$. On utilise à cet effet une m\'ethode de linéarisation locale de l'endomorphisme
\emph{par des homothéties}. 

Techniquement, la difficulté réside dans la 
mise au point de cette m\'ethode de linéarisation car il faut pallier à l'absence du théorème de Koebe.\\
 
Nous présentons maintenant la structure de l'article et les différentes étapes de la démonstration. Les résultats des \emph{sections \ref{lin} et \ref{diam}} concernent la linéarisation et présentent un intérêt pour eux-mêmes. La \emph{section \ref{lin}} est consacrée à la construction d'un procédé de linéarisation général. Il s'agit, pour des choix $\mu$-génériques de $x$, de rendre la suite $(f^n)_n$ normale en $x$ en la 
précomposant par des contractions équivalentes aux applications linéaires
tangentes inverses $(d_xf^{n})^{-1}$. A cet effet, nous estimons
précisément les erreurs cumulées lorsque l'on remplace $f$ par sa différentielle le long 
d'une orbite. Outre la stricte positivité des exposants 
$\lambda_1 \le \cdots \le \lambda_k$, ceci requiert l'hypothèse        
$\lambda_k<2\lambda_1$. Nous obtenons le théorème suivant :

\begin{theo}\label{thli} Si les exposants du système $({\bf P}^k, f, \mu)$ sont tels que 
$\lambda_k < 2 \lambda_1$ alors, pour $\mu$-presque tout point $x$,
la suite $\big(f^n\circ (d_xf^{n})^{-1}\big)_n$ possède au moins une limite injective sur un 
voisinage de $x$.
\end{theo}

En vue d'obtenir un énoncé de linéarisation \emph{par des homothéties}, nous majorons  
la norme des différentielles $(d_x f^{n})^{-1}$. Ceci fait l'objet de la \emph{section \ref{diam}}.
Pour cela, nous reprenons la méthode pluripotentialiste de Briend et Duval \cite{BD} 
dans le contexte des 
 linéarisations. Plus précisément, nous minorons la masse de l'ensemble des points $x$ où les 
 normes
  $\vert\vert(d_xf^{n})^{-1}\vert\vert$ sont ``grandes'' (voir Proposition \ref{THB}). 
  L'énoncé précis de linéarisation suivant résume les informations acquises dans cette section 
  sous une forme
maniable. 

\begin{theo}\label{THL} Si les exposants du système $({\bf P}^k, f, \mu)$ sont tels que 
$\lambda_k < 2 \lambda_1$, alors pour tout borélien $B$, il existe un borélien $\tilde B \subset
B$ de masse arbitrairement proche de $\mu(B)^2$ et $\tau_0>0$ vérifiant les assertions suivantes : pour tout point
$x \in \tilde B$, il existe une suite extraite $(f^{n_j})_j$ et un réel $\nu(x)>0$ tels que 
\begin{itemize}
\item[1.] $f^{n_j}(x) \in B$ pour tout $j\in {\bf N}$.
\item[2.] $f^{n_j} \circ (d_xf^{n_j})^{-1}$ converge uniformément vers un biholomorphisme
sur $B(x,\nu(x))$.
\item[3.] $\vert\vert (d_xf^{n_j})^{-1}\vert\vert \le \tau_0 (\sqrt{d})^{-n_j}$ pour tout $j\in {\bf N}$.
\end{itemize}
\end{theo}

Dans la \emph{section \ref{vol}}, nous montrons que si la mesure $\mu$ est absolument continue
alors les différentielles $(d_xf^{n_j})^{-1}$ intervenant dans le théorème 
\ref{THL} sont équivalentes à des homothéties de rapport $(\sqrt{d})^{-n_j}$. 
Notons que la condition $\lambda_k < 2 \lambda_1$ est 
satisfaite car la régularité de $\mu$ entra\^\i ne la minimalité des exposants.\\

Nous achevons la preuve du théorème \ref{TH} dans la \emph{section \ref{conc}}. 
Nous montrons que  le courant $T$ est régulier en utilisant le procédé de linéarisation 
par les homothéties de rapport $(\sqrt{d})^{-n}$ et les relations d'invariance 
${f^n}^*T=d^n T$. Le résultat de \cite{BL} cité plus haut montre qu'alors $f$ est un 
endomorphisme de Lattès.\\

Nous tenons à remercier
le rapporteur tant pour sa lecture attentive du manuscrit que pour ses conseils de rédaction.

\section{Préliminaires}\label{prel}

Nous résumons ici les principaux outils et résultats utilisés par la suite. 
Nous fixons également quelques notations. 
\subsection{Vocabulaire et notations}

$\bullet$ Un syst\` eme $({\bf P}^k, f, \mu)$ est la donn\'ee d'un endomorphisme holomorphe $f$ de l'espace
projectif de dimension $k$ dont le degr\'e $d$ est sup\'erieur ou \'egal \` a $2$ et dont l'unique mesure
d'entropie maximale est not\'ee $\mu$. Nous dirons aussi que le syst\` eme $({\bf P}^k, f, \mu)$
est de degr\'e $d$.\\

$\bullet$ Soit $(\widehat{{\bf P}^k}, \hat f, \hat\mu)$ l'extension naturelle du système 
$({\bf P}^k, f, \mu)$. On rappelle que  $\widehat{{\bf P}^k}$ est l'ensemble des orbites 
$\left\{\hat x:=(x_n)_{n\in {\bf Z}} \, / \, f(x_n)=x_{n+1} \right \}$ muni de la topologie 
et de la tribu produit. Soient $\pi_0 : \widehat{{\bf P}^k}\to {{\bf P}^k}$ la projection définie par 
$\pi_0(\hat x) = x_0$, $\hat f$ le décalage à droite et  ${\hat f}^{-1}$ le décalage à gauche sur 
$\widehat{{\bf P}^k}$, de sorte que $\pi_0 \circ \hat f =  f \circ \pi_0$. On note $\hat {\mu}$ 
l'unique mesure de probabilité $\hat f$-invariante sur 
$\widehat{{\bf P}^k}$ vérifiant $\pi_{0*}\hat\mu=\mu$. Le caractère mélangeant de $\mu$ passe à 
$\hat \mu$.\\

$\bullet$ Soit $\widehat{X}$ le sous-ensemble de  $\widehat{{\bf P}^k}$ suivant :
\[ \widehat{X}:=\{\hat x \in \widehat{{\bf P}^k} \, / \, x_n \notin Crit(f)  \, , \,  \forall n \in {\bf Z}  \}   \]
où $Crit(f)$ désigne l'ensemble critique de $f$. 
Le borélien $\widehat X$ vérifie $\hat\mu (\widehat X) = 1$, car $\mu$ 
ne charge pas l'ensemble analytique $Crit(f)$ (\cite{S}, Prop. A.6.3).\\ 
Par la suite, on s'autorisera à soustraire à $\widehat X$ 
des ensembles $\hat \mu$-négligeables. \\

\subsection{Branches inverses et exposants}

$\bullet$ On construit une famille de cartes holomorphes $\big(\tau_x\big)_{x\in{\bf P}^k}$ telle que :  
\begin{itemize} 
\item[1.] $\tau_x :{\bf C}^k \to {\bf P}^k$ est un biholomorphisme sur son image et $\tau_x(0)=x$, 
\item[2.] $(\tau_x^*\omega)_0=\frac{i}{2}\sum_{j=1,k} dz_j\wedge d\bar z_j$.
\end{itemize}
où $\omega$ désigne la forme de Fubini-Study. Cette famille est obtenue en explicitant une telle carte en un point base $x_0\in {\bf P}^k$, puis en la propageant à ${\bf P}^k$ par l'action transitive de ${\bf U}_{k+1}({\bf C})$. Ce faisant, on obtient 
plut\^ ot une classe de cartes en $x$ car $\tau_x$ est
définie à un élément du sous-groupe d'isotropie de $x_0$ près.
Cette ambiguïté pourra cependant \^ etre ignorée puisque ${\bf U}_{k+1}({\bf C})$ est compact; les
affirmations faisant intervenir $\tau_x$ devront \^ etre comprises comme valables pour tous les éléments 
de la classe de cartes en $x$. On peut aussi, localement, faire un choix ``différentiable'' 
de $\tau_x$ par rapport à $x$ et en particulier s'assurer que la propriété suivante est vérifiée : 
\begin{itemize}
\item[($\star$)]  
$\tau_{x_0}^{-1}\circ\tau_x - \tau_{x_0}^{-1}(x)$ converge vers l'identité en topologie ${\cal C}^{\infty}$
lorsque $x$ tend vers $x_0$.\\
\end{itemize}

$\bullet$ Nous noterons $B(0,R)$ (resp. $P(0,R)$) la boule euclidienne (resp. le polydisque) de ${\bf C}^k$ centrée en $0$ et de rayon $R$ (resp. de polyrayon $(R,...,R)$). On désignera par $B(x,s)$ l'image de $B(0,s)$ par $\tau_x$. \\

$\bullet$ Nous utiliserons les applications suivantes, où $x\in{\bf P}^k$ et $n \in {\bf N}$:
$$f_x:=\tau_{f(x)}^{-1} \circ f \circ \tau_x$$
$$f_x^n=\tau_{f^n(x)}^{-1} \circ f^n \circ \tau_x=f_{f^n(x)}\circ \cdots \circ f_x.$$   
Elles sont définies sur un voisinage de l'origine de ${\bf C}^k$, dont la taille dépend de $x$ et de $n$. Pour tout $\hat x\in \widehat{X}$, on note  $f^{-n}_{\hat x}$  la branche inverse de $f^n$ ``le long de l'orbite $\hat x$'', 
c'est-à-dire :
$$f^{-n}_{\hat x}:=f^{-1}_{x_{-n}}\circ \cdots  \circ f^{-1}_{x_{-1}}.$$
Le lemme suivant stipule que ces branches inverses existent sur un voisinage de l'origine dont la taille dépend mesurablement de $\hat x$. On trouvera une preuve dans l'article de Briend-Duval \cite{BD} (voir aussi \cite{Du} pp 19-22).
\begin{lem}\label{lempr}
Soient $0<\epsilon \ll 1$ et $0< r_0 \ll 1$. Il existe des fonctions $\rho$ et $r$ continues
sur ${\bf P}^k$ strictement positives hors de $Crit(f)$, ainsi que des fonctions mesurables $\eta : \widehat{X} \to ]0,r_0]$ et $C : \widehat{X} \to [1,+\infty[$  vérifiant les propriétés suivantes :
\begin{itemize} 

\item[1.] Pour tout $x\in {\bf P}^k \setminus Crit(f)$, $f_x$ est injective sur $B(0,\rho(x))$ et 
\[ B(0,r(x)) \subset f_x\big[B(0,\rho(x))\big]. \]
\item[2.] Pour tout $\hat x \in \widehat{X}$, $\lim_n {1 \over n} \log \, \rho(x_n) = 0$.
\item[3.] Pour tout $\hat x\in\widehat{X}$ et tout $n \in {\bf N}$, $f^{-n}_{\hat x}$ est injective sur $B(0,\eta(\hat x))$ et 
$$ \forall \gamma \in ]0,1] \, , \, d_0 f^{-n}_{\hat x}\big[B(0,\gamma\eta (\hat x))\big] \subset B\big(0,\gamma r(x_{-(n+1)})e^{-n(\lambda_1-\epsilon)}\big).$$
\item[4.] $\Lip f^{-n}_{\hat x} \le C(\hat x) e^{-n(\lambda_1-\frac{\epsilon}{2})}$ sur 
$B(0,\eta (\hat x))$. \\
\end{itemize}
\end{lem}

$\bullet$ Les exposants de Liapounov de $\mu$ seront
not\'es $\lambda_1 \le \cdots \le \lambda_k$.  
Nous utiliserons de manière cruciale la minoration optimale de ces exposants :\\

\noindent{\bf Théorème} (Briend-Duval \cite{BD})  
{\it Les exposants d'un système $({\bf P}^k, f, \mu)$ de degré $d$ sont plus grands que $\log \sqrt d$}.

\section{Un proc\'ed\'e de linéarisation}\label{lin}

Notre objectif est de démontrer le théorème \ref{thli} présenté dans l'introduction. Nous adoptons la définition suivante : 

\begin{defi}
Un système $({\bf P}^k, f, \mu)$ est dit linéarisable si pour $\mu$-presque tout
$x \in {\bf P}^k$, il existe $\nu(x)>0$  et une sous-suite de $[f^n \circ \tau_x \circ (d_0f^n_x)^{-1} ]_n$ qui converge uniformément vers une limite injective sur $B(0,\nu(x))$.    
\end{defi}

La proposition suivante fournit deux conditions suffisantes de linéarisibilité. 
La première réduit le problème au contr\^ole uniforme local de la suite 
$f^n_x \circ (d_0f^n_x)^{-1}$ gr\^ace au théorème de Montel. 
La seconde transfère cette question de contr\^ole uniforme en 
"temps négatif", c'est à dire aux applications $f^{n}_{x_{-n}} \circ d_0 f^{-n}_{\hat x}$. 
Nous utilisons pour cela un argument classique basé sur l'invariance de la mesure $\hat\mu$.

\begin{prop}\label{proptech} Soit $({\bf P}^k, f, \mu)$ un système et $R_0 $ un nombre r\'eel strictement
positif. Pour tout $\rho\in ]0,1]$ et $n\in {\bf N}$, on définit les ensembles :
\[  {\cal B}_n(\rho) := \left \{x\in{\bf P}^k \, / \, f^n_x \circ (d_0f^n_x)^{-1} \textrm{ est injective de } B(0,\rho) \textrm{ dans } B(0,R_0) \right\}  \]
\[   {\cal B}(\rho):=\limsup_n{\cal B}_n(\rho).  \]
Le système est linéarisable si l'une des deux conditions suivantes est réalisée :

\noindent 1) Il existe $\alpha : ]0,1] \to {\bf R}^+$ telle que $\lim_{\rho \to 0} \alpha(\rho) = 1$ et $\mu[{\cal B}_n(\rho)]\ge\alpha(\rho)$ pour tout $n\in {\bf N}$.\\
\noindent 2) Pour tout $r_0 \in ]0,R_0]$ il existe des fonctions mesurables $\eta, S : \widehat{X} \to ]0,r_0]$ telles que 
\begin{itemize} 
\item[(i)] $S\le \eta$.
\item[(ii)] Pour tout $\hat x \in \widehat{X}$, $f^{-n}_{\hat x}$ est injective sur $B(0,\eta(\hat x))$.
\item[(iii)] Pour tout $\hat x \in \widehat{X}$ et tout $n\in {\bf N}$,  $d_0f^{-n}_{\hat x}\big[B(0,S(\hat x))\big] \subset
f^{-n}_{\hat x}\big[B(0,\eta(\hat x))\big]$.  
\end{itemize}
La seconde assertion implique la premi\` ere.

\end{prop}

\underline{Démonstration} : 

La linéarisabilité en $x$ résulte, via le théorème de Montel, de l'appartenance de $x$ à 
$\cup_{0<\rho\le 1}{\cal B}(\rho)$. Ainsi, comme $\mu[{\cal B}(\rho)] \ge 
\limsup_n \mu[{\cal B}_n(\rho)]$, la condition 1 entra\^\i ne la linéarisabilité $\mu$-presque partout.\\ 

Voyons maintenant comment la seconde condition entra\^\i ne la première. Posons 
$\widehat{\cal S}(\rho):=\{\hat x \in \widehat{X} / S(\hat x) \ge \rho  \}$. 
Il suffit d'établir les inclusions suivantes :
\begin{equation*}
 \forall n \in {\bf N} \, , \, \pi_0 \big[ \hat{f}^{-n}\big(\widehat{\cal S}(\rho) \big) \big] 
 \subset {\cal B}_n(\rho).
\end{equation*}
En effet, compte tenu de l'invariance de $\hat \mu$, on a 
$\mu\big[{\cal B}_n(\rho) \big]\ge \hat\mu
\big[ \hat{f}^{-n}\big(\widehat{\cal S}(\rho) \big)\big] = 
\hat\mu\big[\widehat{\cal
S}(\rho)\big]$. La fonction $\alpha(\rho):=\hat\mu\big[\widehat{\cal S}(\rho)\big]$ 
convient car $S$ est strictement positive $\hat\mu$-presque partout.

\'Etablissons maintenant les inclusions annoncées. Soit $\hat y := \hat{f}^{n}(\hat x)$
tel que $\hat y \in \widehat{\cal S}(\rho)$. Il s'agit de vérifier que   
$x_0 \in {\cal B}_n(\rho)$. Rappelons que $x_0=\pi_0(\hat x)$. L'appartenance de $\hat y$ à $\widehat{\cal S}(\rho)$ signifie :

$$ d_0f^{-n}_{\hat y}\big[B(0,\rho)\big] \subset d_0f^{-n}_{\hat y}\big[B(0,S(\hat y))\big]
\subset f^{-n}_{\hat y}\big[B(0,\eta(\hat y)\big].$$
Comme $f^{-n}_{\hat y}$ est injective sur $B(0,\eta(\hat y))$ d'inverse $f^n_{x_0}$, on obtient en composant les inclusions précédentes par $f^n_{x_0}$ :
$$f^n_{x_0} \circ \big(d_0 f^n_{x_0} \big)^{-1}\big[B(0,\rho)\big] \subset B(0,\eta(\hat y)) \subset B(0,R_0).$$
Le point $x_0$ appartient donc à ${\cal B}_n(\rho)$. \fin \\

Nous démontrerons le théorème \ref{thli} en vérifiant que la condition 2 de la proposition \ref{proptech}
est satisfaite. Ceci consistera à compenser les erreurs dues à la substitution de        
$d_0 f^{-1}_{x_j}$ à $ f^{-1}_{x_j}$ le long de $\hat x$ en diminuant le rayon $\eta(\hat x)$. Pour que
les compensations cumulées fournissent un rayon $S(\hat x)$ strictement positif, les erreurs commises devront \^etre
négligeables devant la plus petite dimension caractéristique de l'ellipsoïde $d_0 f^{-j}_{\hat
 x}\big[B(0,1)\big]$. L'objet  du lemme suivant est de montrer que tel est le cas lorsque les exposants vérifient l'inégalité  
 $\lambda_k < 2\lambda_1$.

\begin{lem}\label{lem3.4}
Soient un système $({\bf P}^k, f, \mu)$ et $0 < \epsilon \ll 1$. Il existe des fonctions mesurables $\eta, E, F : 
 \widehat{X} \to ]0,+\infty[$ vérifiant $0< \eta \le r_0 \le R_0$ telles que pour tout 
 $\hat x=(x_n)_{n\in {\bf Z}}$ \'el\'ement de $\widehat{X}$ et tout $n\in {\bf N}$ : \\
\noindent 1)  $f^{-n}_{\hat x}$ est injective sur $B(0,\eta(\hat x))$.\\    
\noindent 2) Pour tout  $\gamma\in ]0,1]$ et tout $u \in 
d_0f^{-n}_{\hat x}\big[B(0,\gamma\eta(\hat x))\big]$ : 
$$ \vert\vert \big(d_0f^{-1}_{ x_{-(n+1)}} - f^{-1}_{ x_{-(n+1)}}\big) (u) \vert\vert \le \gamma E(\hat x) 
e^{-2n(\lambda_1-\epsilon)}.$$
\noindent 3)  $\vert\vert d_0f^{n+1}_{ x_{-(n+1)}}\vert\vert \le F(\hat x)e^{n(\lambda_k+\epsilon)}.$          
\end{lem}

\underline{Démonstration} :

Nous utilisons ici le lemme \ref{lempr}. Pour tout  $\hat x \in \widehat{X}$, l'assertion 1 est satisfaite. De plus, l'application $f_{ x_{-(n+1)}}$ est inversible sur $B(0,r)$, où $r:=r(x_{-(n+1)})$, et son inverse $g$ est à valeurs dans $B(0,\rho)$ où  $\rho:=\rho\big(x_{-(n+1)} \big)$. Soit $\sum_{p\ge 2} Q_p$ le développement de Taylor de $g-d_0g$, où $Q_p$ désigne une application 
homogène de degré $p$. Si $u\in B(0,r)$ alors $\vert\vert Q_p(u)\vert\vert = \vert\vert\frac{1}{2\pi}
\int_{0}^{2\pi}g(e^{i\theta}u)e^{-ip\theta}d\theta \vert\vert \le \rho$ et donc :

$$ \vert\vert (g-d_0 g) (u)\vert\vert \le \sum_{p\ge 2} \frac{\vert\vert u \vert\vert ^p}{r^p}
\vert\vert Q_p\big( \frac{ru}{\vert\vert u \vert\vert}\big)\vert\vert \le
 \rho \sum _{p\ge 2} \left( \frac{\vert\vert u \vert\vert }{r} \right)^p.$$
Lorsque de plus $u \in d_0f^{-n}_{\hat x}\big[B(0,\gamma\eta({\hat x})\big]$ alors $\frac{\vert\vert u \vert\vert}{r} \le \gamma e^{-n(\lambda_1 - \epsilon)}$ (cf lemme \ref{lempr},(3)) et il s'ensuit que :
$$\vert\vert (g-d_0 g) (u)\vert\vert \le \frac{\gamma^2\rho}{1-\gamma e^{-(\lambda_1-\epsilon)}} e^{-2n(\lambda_1-\epsilon)}.$$   

L'assertion 2 du lemme s'en déduit car $\rho(x_{-(n+1)})$ a un taux de croissance exponentiel nul (cf lemme \ref{lempr},(2)). La dernière assertion découle immédiatement de la définition des exposants de Liapounov. Nous \^otons ici à $\widehat X$ un sous-ensemble de $\hat\mu$ mesure nulle. \fin \\

\underline{Démonstration du théorème \ref{thli}} : \\
Il s'agit de montrer que la condition 2 de la proposition \ref{proptech} est satisfaite
lorsque $\lambda_k<2\lambda_1$. Reprenons les notations du lemme \ref{lem3.4} et introduisons sur $\widehat X$ les fonctions 
mesurables suivantes :

$$\xi_n(\hat x):= \sup \left\{  t\le \eta(\hat x) \, / \, d_0f^{-n}_{\hat x}\big[B(0,t)\big] \subset
f^{-n}_{\hat x}\big[B(0,\eta({\hat x}))\big]  \right\}$$ 

$$n_0(\hat x) := \min  \left \{p \ge 1  \, /  \, \forall n \ge p  \, :  \, \frac{EF}{\eta}(\hat x)\le e^{n\epsilon}  \right  \}$$

$$s(\hat x) := \min   \left \{\xi_n(\hat x) \, : \, 0\le n \le n_0(\hat x) \right\}.$$
Posons $\kappa_j:=1-e^{-j(2\lambda_1-\lambda_k-6\epsilon)}$ avec $\epsilon$ suffisamment petit pour que le produit $\prod_{j=1}^{\infty} \kappa_j$ converge et soit strictement positif. Définissons les fonctions $s_n$ par :
\[ s_n(\hat x) := s(\hat x)   \textrm{ si }  n\le n_0(\hat x) \]
\[ s_n(\hat x) := s(\hat x)\prod_{j=n_0(\hat x)}^{n-1}\kappa_j  \textrm{ si }  n\ge n_0(\hat x)+1. \]
Pour montrer que la fonction $S(\hat x):=s(\hat x)\prod_{j=1}^{\infty}\kappa_j$ convient, il suffit d'établir les inclusions :  
$$ (I_n)_{n\ge 0}\;\;\;:\;\;\;d_0f^{-n}_{\hat x}\big[B(0,s_n(\hat x))\big]
\subset
f^{-n}_{\hat x}\big[B(0,\eta({\hat x}))\big].$$   
Par définition de $s_n(\hat x)$, ces inclusions sont satisfaites lorsque $n\le n_0(\hat x)$. Supposons que $(I_n)$ soit vraie pour $n\ge n_0(\hat x)$ et posons $\nu_n:= \big(\frac{E s_n}{\eta}\big)(\hat x)e^{-2n(\lambda_1-2\epsilon)}$. On a alors : 
\begin{equation}\label{E1}
 s_{n+1}\le s_n - \vert\vert\big( d_0f^{-(n+1)}_{ \hat x}\big)^{-1}\vert\vert \nu_n.
\end{equation}
En effet :
$$s_{n+1}=s_n\kappa_n=s_n\big(1-e^{-n(2\lambda_1-\lambda_k -6\epsilon)}\big)\le
s_n\big(1-\frac{EF}{\eta}(\hat x)e^{-n(2\lambda_1-\lambda_k -5\epsilon)}\big)$$
$$\le
s_n - \vert\vert d_0f^{n+1}_{ x_{-(n+1)}}\vert\vert\nu_n = 
s_n - \vert\vert\big( d_0f^{-(n+1)}_{ \hat x}\big)^{-1}\vert\vert \nu_n$$ 
la première majoration résultant de la définition de $n_0(\hat x)$ et la seconde du lemme \ref{lem3.4},(2). 

Désignons par $\Lambda$ la frontière de $d_0f^{-(n+1)}_{\hat x}\big[B(0,s_n)\big]$. On vérifie aisément que l'inégalité (\ref{E1}) se traduit par :
\begin{equation}\label{E2}
d_0f^{-(n+1)}_{\hat x}\big[B(0,s_{n+1})\big] \subset
d_0f^{-(n+1)}_{\hat x}\big[B(0,s_n)\big] \backslash \bigcup_{p\in \Lambda}B(p,\nu_n).
\end{equation}
Par ailleurs, la première assertion du lemme \ref{lem3.4} (où l'on prend $\gamma=\frac{s_n}{\eta}$)
stipule que sur $d_0f^{-n}_{\hat x}\big[B(0,s_{n})\big]$, $f^{-1}_{ x_{-(n+1)}}$ diffère d'au plus $\nu_n$ de sa différentielle. Il s'ensuit que 
\begin{equation}\label{E3}
d_0f^{-(n+1)}_{\hat x}\big[B(0,s_n)\big] \backslash \bigcup_{p\in \Lambda}B(p,\nu_n)
\subset f^{-1}_{ x_{-(n+1)}} \circ d_0f^{-n}_{\hat x}\big[B(0,s_{n})\big].  
\end{equation}
Observons finalement que l'inclusion $(I_n)$, composée par $f^{-1}_{ x_{-(n+1)}}$, s'écrit
\begin{equation}\label{E4}
f^{-1}_{ x_{-(n+1)}} \circ d_0f^{-n}_{\hat x}\big[B(0,s_{n})\big]
\subset
f^{-(n+1)}_{\hat x}\big[B(0,\eta({\hat x})\big].
\end{equation}
Les inclusions (\ref{E2}), (\ref{E3}) et (\ref{E4}) encha\^\i nées donnent $(I_{n+1})$. \fin

\begin{rem}
L'inégalité  $\lambda_k < 2\lambda_1$ peut \^etre interpr\^etée comme une condition de \emph{non résonance} 
entrainant la linéarisabilité. Jonsson et 
Varolin (cf \cite{JV}, Thm 3) ont, indépendamment de nous, mis en évidence la m\^eme condition dans un
problème voisin.
\end{rem}

\section{Une version précisée du proc\'ed\'e de linéarisation}\label{diam}

L'objet de cette section est de contr\^oler le diamètre des ellipsoïdes $(d_0f^n_x)^{-1}\big[B(0,1)\big]$ associés au procédé de linéarisation fourni par le théorème \ref{thli}.  Nous en déduisons le théorème \ref{THL} énoncé dans l'introduction. 

Le théorème de Briend-Duval, déjà utilisé implicitement  pour établir le lemme \ref{lempr}, majore 
le taux de décroissance exponentielle de la taille 
de ces ellipsoïdes par $-\log \sqrt d$. Cela signifie que pour tout $\epsilon >0$, on a
$\vert\vert (d_0f^n_x)^{-1} \vert\vert \lesssim e^{n\epsilon} (\sqrt d)^{-n}$ pour $n$ assez grand.
En reprenant la méthode de Briend-Duval dans le contexte de la proposition \ref{proptech}, nous 
obtenons une majoration plus précise :  
$\vert\vert (d_0f^n_x)^{-1} \vert\vert \lesssim  (\sqrt d)^{-n}$. Rappelons que  ${\cal B}_n(\rho)$ 
est défini par :
\[  {\cal B}_n(\rho) := \left \{x\in{\bf P}^k \, / \, f^n_x \circ (d_0f^n_x)^{-1} \textrm{ est 
injective de } B(0,\rho) \textrm{ dans } B(0,R_0) \right\}  \]
et qu'en vertu de la proposition \ref{proptech} et de la preuve du th\'eor\` eme \ref{thli}, il existe une fonction
$\alpha: ]0,1] \to {\bf R}^+$ telle que $\lim_{\rho \to 0} \alpha(\rho) = 1$ et 
\[\mu[{\cal B}_n(\rho)] \ge \alpha(\rho).\]  
Nous montrons la proposition suivante :
\begin{prop}\label{THB} Soit $({\bf P}^k, f, \mu)$ un système de degré $d \geq 2$ tel que  
$\lambda_k < 2 \lambda_1$. On pose pour  $\tau >0$, $\rho \in ]0,1]$ et $n\in{\bf N}$ : 
\[  {\cal D}_n(\rho,\tau) :={\cal B}_n(\rho)\cap  \left\{ x\in{\bf P}^k \, / \, 
\vert\vert (d_0f^n_x)^{-1} \vert\vert  \le \tau (\sqrt d)^{-n} \right \} .\]
Alors on a l'inégalité :
\[ \liminf_n \mu[{\cal D}_n(\rho,\tau)] \ge \alpha(\rho)-\frac{C}{\tau^2\rho^2}. \]
où $C>0$ et $\alpha: ]0,1] \to {\bf R}^+$ est une fonction telle que  
$\lim_{\rho \to 0} \alpha(\rho) = 1$. 

\end{prop}

Le principe de la preuve est le suivant. Puisque $\mu[{\cal B}_n(\rho)] \ge \alpha(\rho)$ d'après la proposition \ref{proptech}, il s'agit de majorer la mesure du complémentaire de ${\cal D}_n(\rho,\tau)$ dans ${\cal B}_n(\rho)$, noté ${\cal D}^c_n(\rho,\tau)$. Or, par tout point de ${\cal D}^c_n(\rho,\tau)$ passe un disque dont le diamètre est 
au moins égal à $\tau \rho (\sqrt d)^{-n}$ et dont l'image par $f^n$ reste contenue dans une boule
de rayon $R_0$ fixé. Comme $f^{n*}T=d^n T$, il passe donc par tout point de ${\cal D}^c_n(\rho,\tau)$ un ``grand'' disque ``peu'' chargé par $T$. Des techniques pluripotentialistes permettent alors de majorer précisément la masse de l'ensemble de ces points pour la mesure $\mu=T^k$.\\

\underline{Démonstration de la proposition \ref{THB}} :

On dira qu'un disque holomorphe $\sigma : \Delta \to {\bf C}^k$ est de taille $l > 0$ et passe par 
$z\in {\bf C}^k$ si il est 
de la forme $\sigma(u)=z+l u. {v} +\beta(u)$ où ${v}$ est un vecteur unitaire de ${\bf C}^k$, 
$\beta(0)=0$ et $\vert\vert \beta \vert\vert\le \frac{l}{1000}$. \\

L'ingrédient principal est le théorème suivant dont la preuve est résumée dans l'appendice : \\

\noindent{\bf Théorème} (Briend-Duval \cite{BD}) 
{\it Soit $S:=dd^c w$ un $(1,1)$-courant positif fermé de potentiel $w$ continu sur le
polydisque $P(0,R)$ et 
$E\subset P(0,\frac{R}{2})$. On suppose que par tout point $z \in E$ passe un disque holomorphe 
$\sigma_z : \Delta
\to {\bf C}^k$ de taille $l$ et qu'il existe une fonction $h_z$ harmonique sur $\Delta$
telle que 
$\vert w \circ \sigma_z - h_z \vert \le \epsilon$ sur $\Delta$. Alors il existe une constante $C(w)$ ne
dépendant que de $w$ telle que $S^k(E) \le C(w) \frac{k^2}{l^2}\epsilon$}.\\

En vue d'utiliser ce résultat, nous fixons des systèmes de coordonnées locales sur ${\bf P}^k$.
Considérons un recouvrement de ${\bf P}^k$ par des ouverts $U_1,...,U_N$ centrés en des points $m_j$
et tel que sur chaque $U_j$ nous puissions fixer des déterminations des cartes $\tau_x$ dépendant
différentiablement de $x$ (cf la condition $(\star)$, section \ref{prel}).
Posons $\tau_j:=\tau_{m_j}$ puis, pour $R>0$ fixé, $V_j:=\tau_j(P(0,R))$. 
Si le recouvrement est assez
fin alors les propriétés suivantes sont satisfaites :

\begin{itemize} 

\item[(i)] $U_j \subset \tau_j(P(0,\frac{R}{2}))$ et $\tau_x(P(0,\frac{R}{2}))\subset V_j$ pour tout $x\in U_j$
\item[(ii)] $\forall x \in U_j, 
\vert\vert \tau_j^{-1}\circ \tau_x - (\tau^{-1}_j(x)+\Id) 
\vert\vert_{{\cal C}^1,\overline{P(0,\frac{R}{2})}}
\le\frac{1}{1000}$

\end{itemize}
puis, si $R_0$ (introduit au lemme \ref{lempr}) est pris assez petit :

\begin{itemize}
\item[(iii)] 
$\forall x \in {\bf P}^k, \exists l \in \{1,...,N\}$ tel que $\tau_x\big[B(0,R_0)\big]
\subset V_l$  
\item[(iv)] $\mu \{x \in U_j \cap {\cal B}_n(\rho)\;/\;     
(d_0f^n_x)^{-1}\big[B(0,\rho)\big] \subset P(0,\frac{R}{2})\}=
\mu (U_j \cap {\cal B}_n(\rho)) - \epsilon_{n,j}$ avec 
$\lim_n \epsilon_{n,j}=0$
\end{itemize}
enfin, si $v_j$ désigne un potentiel continu de $T$ sur $V_j$, il existe une constante $M>0$ telle que :  

\begin{itemize}
\item[(v)]  $T = dd^c v_j$ et $\vert v_j\vert \le M$ sur $V_j$ pour tout $j \in \{1,...,N\}$. 
\end{itemize}

D'après le th\'eor\` eme \ref{thli}, il existe une fonction $\alpha$ qui v\'erifie la propri\'et\'e
\'enonc\'ee \` a la proposition \ref{proptech} (1) ; autrement dit, $\alpha(\rho)$ tend vers $1$ quand 
$\rho$ tend
vers $0$ et $\mu[{\cal B}_n(\rho)] \ge \alpha(\rho)$.\\    
Comme il s'agit 
de minorer
$\liminf_n\mu[{\cal D}_n(\rho,\tau)]$, la propriété (iv) montre que l'on peut considérer que :
\begin{equation}\label{inclu}
 \forall x \in U_j \cap {\cal B}_n(\rho) \ , \ (d_0f^n_x)^{-1}\big[B(0,\rho)\big] \subset P(0,\frac{R}{2}).
\end{equation}
Pour tout $j \in \{1,...,N\}$
nous allons établir que : 
\begin{equation}\label{majo} \mu\big[{\cal D}^c_n(\rho,\tau)\cap U_j\big] \le
C(v_j\circ\tau_j)\frac{Mk^2}{\tau^2\rho^2}.
\end{equation}
Soit donc 
$x\in{\cal D}^c_n(\rho,\tau)\cap U_j$ et $V_n(x)$ un vecteur unitaire tel que 
$\vert\vert (d_0f^n_x)^{-1} \vert\vert=\vert\vert (d_0f^n_x)^{-1} [V_n(x)] \vert\vert$. 
Si on note $v_n(x)=(d_0f^n_x)^{-1} [V_n(x)]$, on a 
$\vert\vert v_n(x)\vert\vert \ge \tau \rho (\sqrt{d})^{-n}$. On définit ainsi un disque affine 
$\Phi_{n,x}: \bar{\Delta}\to{\bf C}^k$
de diamètre au moins égal à $\tau \rho  (\sqrt{d})^{-n}$ en posant : 
$$\Phi_{n,x}(t):=(d_0f^n_x)^{-1}[t\rho. V_n(x)]=t\rho. v_n(x).$$ Comme
$x\in U_j$, (\ref{inclu}) et (i) permettent
de définir un nouveau disque $\Phi_{j,n,x}:\Delta \to P(0,R)$ par     
$\Phi_{j,n,x}:=\tau_j^{-1}\circ \tau_x \circ \Phi_{n,x}$. Compte tenu de la propriété (ii),   
$\Phi_{j,n,x}$ est un disque holomorphe de taille $l:=\rho\vert\vert v_n(x) \vert\vert 
\ge \tau \rho  (\sqrt{d})^{-n}$
passant par $\tau_j^{-1}(x)$.\\
Choisissons $l \in \{1,...,N\}$ tel que $\tau_{f^n(x)}\big[B(0,R_0)\big]
\subset V_l$ (propriété (iii)) alors, comme $x \in {\cal B}_n(\rho)$, on a $f^n \circ \tau_x \circ 
\Phi_{n,x}(\Delta) \subset \tau_{f^n(x)}\big[B(0,R_0)\big]
\subset V_l$ et donc $$dd^c(v_l\circ f^n \circ \tau_x \circ 
\Phi_{n,x})= (f^n \circ \tau_x \circ 
\Phi_{n,x})^*T=(\tau_x \circ 
\Phi_{n,x})^*f^{n*}T=d^n(\tau_x \circ 
\Phi_{n,x})^*T.$$
 Par ailleurs, puisque $\tau_x\circ \Phi_{n,x}(\Delta)\subset V_j$ (cf. (\ref{inclu}) et (i)), on a         
$$(\tau_x \circ 
\Phi_{n,x})^*T=dd^c(v_j\circ \tau_x \circ 
\Phi_{n,x})=dd^c(v_j\circ \tau_j \circ 
\Phi_{j,n,x}).$$
Ainsi, $dd^c\big[v_l\circ f^n \circ \tau_x \circ 
\Phi_{n,x} - d^n v_j\circ \tau_j \circ 
\Phi_{j,n,x}\big] =0$.
Autrement dit, la fonction entre crochets est harmonique sur $\Delta$ et, 
puisque $\vert v_l\vert\le M$, le
potentiel $v_j\circ\tau_j$ de $\tau_j^* T$ diffère d'au plus $\frac{M}{d^n}$ d'une fonction
harmonique $h$ sur le disque 
$\Phi_{j,n,x} $ de taille $l\ge \tau \rho  (\sqrt{d})^{-n}$. 
Dans ces conditions, (\ref{majo}) découle
immédiatement
du théorème de Briend-Duval. On en déduit l'estimation annoncée avec 
$C=Mk^2\sum_{j = 1}^{N}C(v_j\circ\tau_j)$.  \fin \\

Terminons cette section par la preuve du théorème \ref{THL}. Il s'agit d'établir une version du procédé de linéarisation où les orbites issues d'un borélien prescrit sont assujetties à récurrence. Cette précision découle des estimations fournies par la proposition \ref{THB} et du caractère mélangeant de $\mu$. \\

\underline{Démonstration du théorème \ref{THL}} :

\noindent Posons ${\cal D}_n(\rho,\tau,B):={\cal D}_n(\rho,\tau)\cap B \cap f^{-n}(B)$
et ${\cal D}(\rho,\tau,B):=\limsup_n {\cal D}_n(\rho,\tau,B)$. Il est clair que si 
$x\in{\cal D}(\rho_0,\tau_0,B)$ alors il existe une suite extraite $(f^{n_j})_j$ vérifiant
les trois assertions du théorème \ref{THL}. Il suffit donc     
d'observer que $\mu \big({\cal D}_n(\rho,\tau,B) \big)$ approche $\mu(B)^2$ par défaut
pourvu que $\rho_0$,$\frac{1}{\tau_0}$ soient assez petits et $n$ assez grand.
Or ceci résulte immédiatement de la proposition \ref{THB} et du caractère mélangeant de $\mu$.
Il suffit en effet de fixer $\rho_0$ assez petit puis $\tau_0$ assez grand pour que 
${\cal D}_n(\rho,\tau)$ soit presque de $\mu$-mesure pleine pour $n$ assez grand et d'utiliser ensuite le
fait que $\mu\big[B \cap f^{-n}(B)\big]$ approche $\mu(B)^2$ lorsque $n$ tend vers l'infini.
   \fin \\

\section{Linéarisation par des homothéties}\label{vol}

Dans cette partie, nous montrons que la suite des itérées $(f^n)_n$ est linéarisable par
 des homothéties de
rapport $(\sqrt{d})^{-n}$ si et seulement si $\mu$ est absolument continue par rapport
 à la mesure de Lebesgue.  Nous adoptons la définition suivante :
\begin{defi}\label{defi4.1}
Un système $({\bf P}^k, f, \mu)$ de degré $d$ est dit $\sqrt{d}$-linéarisable si
 pour $\mu$-presque tout
$x \in {\bf P}^k$, il existe $\nu(x)>0$  et une sous-suite de 
$[f^n \circ \tau_x \circ  (\sqrt{d})^{-n} \Id ]_n$ qui converge uniformément vers une limite injective
 sur $B(0,\nu(x))$.    
\end{defi} 

Autrement dit, un système est $\sqrt{d}$-linéarisable si pour tout $x$ générique, les ellipsoïdes 
$(d_0f^{n}_x)^{-1}\big[B(0,1)\big]$ sont assimilables à des boules euclidiennes de rayon $(\sqrt{d})^{-n}$. Comme la taille de ces ellipsoïdes est au plus $\tau_0 (\sqrt{d})^{-n}$ (cf  théorème \ref{THL}), il suffit d'en contr\^oler le volume. L'absolue continuité de $\mu$ le permet. Nous introduisons à cet effet les ensembles suivants  :
\[  \forall \nu \in ]0,1] \ , \ {\cal V}_n(\nu):= \left\{  x\in {\bf P}^k  \, / \, 
\nu^2 d^{kn} \le \vert \Jac f^n_x\vert^2 \le \frac{1}{\nu^2} d^{kn} \right\}  \]
où $\Jac f^n_x$ désigne le Jacobien complexe de $f^n_x$ en $0$. Nous obtenons le résultat suivant :
\begin{prop}\label{pr4.1}
Soit $({\bf P}^k, f, \mu)$ un système de degré $d$. Les propriétés suivantes sont équivalentes :
\begin{itemize} 
\item[1.] $\mu$ est absolument continue par rapport à la mesure de Lebesgue $\omega^k$.
\item[2.] Les exposants du système sont tous égaux à $\log \sqrt{d}$  et il existe 
$\beta : ]0,1] \to {\bf R}^+$ vérifiant  $\lim_{\nu \to 0} \beta(\nu) = 1$ et 
$ \liminf_n\mu[{\cal V}_n(\nu)]\ge\beta(\nu)$.
\item[3.] Le système est $\sqrt{d}$-linéarisable.
\end{itemize}
\end{prop}

Nous noterons $Y$ l'ensemble $\cup_p Crit\;f^p$. En tant qu'union d\'enombrable de sous-vari\'et\'es
alg\'ebriques de ${\bf P}^k$, c'est un ensemble de $\mu$-mesure nulle (voir \cite{S}) : $\mu(Y)=0$.  

\underline{Démonstration} :

$1) \Rightarrow 2)$. Commençons par établir l'existence de la fonction $\beta$. On note 
$m=\omega^k$ la mesure de Lebesgue sur ${\bf P}^k$. Puisque $\mu$ est absolument continue par rapport \` a
$m$, il existe 
$\varphi \in L^1(m)$ telle que 
$\mu=\varphi \, dm$. D'après le théorème de Lusin, il existe pour tout $n\in {\bf N}$ des
fonctions continues $g_n$ et $h_n$ ainsi que des boréliens $C_n(\varphi)$ et $C_n(\varphi \circ f^n)$ 
vérifiant : 
\[ \varphi = g_n \textrm{ sur }  C_n(\varphi) \ \  \textrm{ et } \ \ \mu\big[C_n(\varphi)\big] 
\ge 1-\frac{1}{n} \]   
\[ \varphi \circ f^n = h_n  \textrm{ sur } C_n(\varphi \circ f^n) \ \  
\textrm{ et } \ \ \mu\big[C_n(\varphi \circ f^n)\big] \ge 1-\frac{1}{n}.  \]
Soit $A_{\nu}:=\{x \in {\bf P}^k\; /\; \nu <\varphi(x) <\frac{1}{\nu}\}$ où 
$\nu \in ]0,1]$. On pose :  
\[ Z_{n,\nu}:=\big[f^{-n}(A_{\nu})\cap A_{\nu}\big] \cap \big[C_n(\varphi) \cap C_n(\varphi \circ f^n)\big]  \cap Y^{c}. \]  
Rappelons que l'ensemble $Z_{n,\nu}^{Leb}$ des points de Lebesgue de $Z_{n,\nu}$ 
est défini par :
$$Z_{n,\nu}^{Leb}:=\{x\in Z_{n,\nu}\;/\;\lim_{s\to 0}
\frac{m\big[B(x,s)\cap Z_{n,\nu}\big]}{m\big[B(x,s)\big]}=1 \} . $$
L'absolue continuité de $\mu$ entra\^\i ne  $\mu\big(Z_{n,\nu}^{Leb}\big)=\mu\big(Z_{n,\nu}\big)$. 
Compte tenu du caractère
mélangeant de $\mu$ et du fait que $\mu(Y)=0$, on obtient pour $n$ assez grand :
$$\mu \big(Z_{n,\nu}^{Leb}\big) \ge \mu(A_{\nu})^2\left(1-{\nu \over 2}\right)  - \frac{2}{n} 
\geq \mu(A_{\nu})^2(1-\nu). $$
La fonction $\beta(\nu):=\mu(A_{\nu})^2(1-\nu)$ convient si l'inclusion 
$Z_{n,\nu}^{Leb}\subset {\cal V}_n(\nu)$ est satisfaite. Fixons donc $x \in Z_{n,\nu}^{Leb}$. 
Puisque $x$ n'appartient pas à $Crit\;f^n$, il existe $s_0>0$ tel que $f^n$ soit injective sur $B(x,s_0)$. En outre, $x$
étant un point de Lebesgue de $Z_{n,\nu}$, on peut diminuer $s_0$ pour que       
$m\big[B(x,s)\cap Z_{n,\nu}\big] \ge \frac{1}{ 2}m\big[B(x,s)\big]>0$ pour tout $0<s<s_0$. En utilisant 
des
changements de
variables, d'abord par rapport à $\mu=\varphi \, dm$ qui est de Jacobien constant égal à $d^k$, puis 
par rapport à
$m=\omega^k$, on obtient :
$$d^{kn}\int_{B(x,s)\cap Z_{n,\nu}}\varphi \, \omega^k = \int_
{f^n[B(x,s)\cap Z_{n,\nu}]}\varphi \, \omega^k
=\int_{B(x,s)\cap Z_{n,\nu}}\varphi \circ f^n \,  (f^{n*}\omega^k).$$
Or, puisque $C_n(\varphi) \cap 
C_n(\varphi \circ f^n)$ contient $Z_{n,\nu}$, on peut remplacer $\varphi$ par $g_n$ et 
$\varphi \circ f^n$ par $h_n$ dans ces
intégrales. Après normalisation par $m(s,n,\nu):= m[B(x,s)\cap Z_{n,\nu}]$, il vient :
$$\frac{d^{kn}}{m(s,n,\nu)}\int_{B(x,s)\cap Z_{n,\nu}}g_n \, \omega^k  =
\frac{1}{m(s,n,\nu)} 
\int_{B(x,s)\cap Z_{n,\nu}}h_n \, (f^{n*}\omega^k).$$ 
 Comme les fonctions $g_n$ et $h_n$ sont continues en $x$ et 
$(f^{n*}\omega^k)_x=\vert \Jac f^n_x\vert^2 (\omega^k)_x$, on obtient lorsque $s$ tend vers $0$ :
$$d^{kn} \varphi(x) = d^{kn} g_n(x)=
h_n(x)\vert \Jac f^n_x\vert^2 = \varphi \circ f^n(x)\vert \Jac f^n_x\vert^2 $$
c'est à dire $\frac{\vert \Jac f^n_x\vert^2}{d^{kn}}=
\frac{\varphi(x)}{\varphi \circ f^n(x)}$. Ainsi $x \in {\cal V}_n(\nu)$ car $x$ et 
$f^n(x)$ appartiennent à $A_{\nu}$.\\

Vérifions maintenant que les exposants de $\mu$ sont minimaux. On dispose de l'égalité classique 
$\lim_n \frac{1}{n}\log \vert \Jac f^n_x\vert^2 =2\sum_{i=1}^{k} \lambda_i$, valable pour 
$\mu$-presque tout $x$ (cf par exemple \cite{A} Section 3.3). Notons 
${\cal V}(\nu):=\limsup_n {\cal V}_n(\nu)$ et choisissons $\nu$ assez petit pour que 
$\mu[{\cal V}(\nu)] \ge \beta(\nu) \ge \frac{1}{2}$. Comme 
$\lim_n \frac{1}{n} \log \vert \Jac f^n_x\vert^2 = k \log d$ pour $x \in {\cal V}(\nu)$, on obtient 
$\sum_{i=1}^{k} \lambda_i = k  \log \sqrt{d}$. La minimalité des exposants découle alors de la minoration 
$\lambda_i \ge  \log \sqrt{d}$.\\

$2) \Rightarrow 3)$. La proposition \ref{THB} s'applique car les exposants sont tous égaux à 
$\log \sqrt d$. Nous en reprenons les
notations et posons 
\[ {\cal DV}_n(\rho,\tau,\nu):= {\cal D}_n(\rho,\tau)\cap{\cal V}_n(\nu) \ \textrm{ et }  
\ {\cal DV}(\rho,\tau,\nu):= \limsup_n {\cal DV}_n(\rho,\tau,\nu).  \]
D'après 2) et la proposition \ref{THB}, 
$\mu[{\cal DV}(\rho,\tau,\nu)]$ est arbitrairement proche de $1$ pourvu que $\rho$ et $\nu$ soient assez
petits et $\tau$ suffisamment grand. Il suffit donc de montrer que $(f_x ^n)_n$ est linéarisable par 
$\Lambda_n:=
(\sqrt{d})^{-n} \Id$ lorsque $x \in {\cal DV}(\rho,\tau,\nu)$. Soit donc $(n_j)_j$ une suite strictement
croissante d'entiers telle que $x\in{\cal DV}_{n_j}(\rho,\tau,\nu)$ pour tout $j$. Puisque     
${\cal DV}_{n_j}(\rho,\tau,\nu) \subset  {\cal D}_{n_j}(\rho,\tau) \subset {\cal B}_{n_j}(\rho)$ on a 
\[ f^{n_j}_x \circ (d_0f^{n_j}_x)^{-1} \big(B(0,\rho)\big) \subset B(0,R_0). \]
Il s'agit donc de vérifier que $(d_0f^{n_j}_x)^{-1}$ est équivalente à $\Lambda_{n_j}$. 
A cet effet, notons 
\[ \delta_{j,1}\le \cdots \le \delta_{j,k} \]
les valeurs singulières de $P := (d_0f^{n_j}_x)^{-1}$, c'est à dire les valeurs propres de 
la racine carrée de $PP^*$, où $P^*$ désigne l'adjoint de $P$.  Il existe en particulier deux 
matrices unitaires $U$ et $V$ telles que $U P V = Diag(\delta_{j,1}, \cdots , \delta_{j,k})$. 
Ces valeurs singulières vérifient $\delta_{j,k}\le \tau (\sqrt{d})^{-n_j}$ car 
$x \in {\cal D}_{n_j}(\rho,\tau)$
et $(\delta_{j,1}...\delta_{j,k})^2 = \vert \Jac f^{n_j}_x \vert^{-2} \ge \nu^2 d^{-kn_j}$ car
$x\in {\cal V}_{n_j}(\nu)$. D'où l'on déduit les inégalités :
\[ \nu \tau^{1-k} (\sqrt{d})^{-n_j} \le\delta_{j,1}\le ...\le \delta_{j,k}\le \tau (\sqrt{d})^{-n_j}. \] 
L'application $(d_0f^{n_j}_x)^{-1}$ est donc équivalente à l'homothétie $\Lambda_{n_j}$. \\ 

$3) \Rightarrow 1)$. Soit $x\in {\bf P}^k$ un point $\mu$ générique. D'après 3), il existe $\rho >0$ et 
une suite croissante d'entiers  $(n_j)_j$ tels que $f^{n_j} \circ \tau_x \circ \Lambda_{n_j} : B(0,\rho) 
\to B(0,R_0)$ soit une suite d'injections. Soient $B_r := \tau_x [ B(0,r)]$ et $B_{n_j}:= 
\tau_x [ B(0,\rho (\sqrt{d})^{-n_j})]$. Il s'ensuit que :
$$\liminf_{r\to 0} \frac   {\mu (B_r)}   {m (B_r)}   \le
\liminf_j          \frac   {\mu (B_{n_j})}  {m (B_{n_j})}  \lesssim 
\liminf_j          \frac   {\mu (B_{n_j})}  {d^{-kn_j}} =\liminf_j \mu (f^{n_j} (B_{n_j})) \le 1$$
où la dernière égalité provient du fait que $\mu$ est de jacobien constant $d^k$. Ceci étant vérifié pour 
$\mu$-presque tout $x$, la mesure $\mu$ est bien absolument continue par rapport à $m$ (cf \cite{Ma}, Th 2.12). \fin    

\begin{rem}\label{remrec}
Comme nous l'avons fait pour établir le théorème \ref{THL}, une légère modification dans la preuve de $2) 
\Rightarrow 3)$ permet de choisir la sous-suite $(n_j)_j$ de fa\c con à ce que $f^{n_j}(x)$ ne s'échappe 
pas d'un borélien $B$ de $\mu$-mesure strictement positive prescrit.  
\end{rem}

\section{Régularisation du courant de Green} \label{conc}

Nous achevons ici la preuve du théorème \ref{TH}. D'après la proposition \ref{pr4.1},
il s'agit de caractériser les systèmes $({\bf P}^k, f, \mu)$ qui sont    
$\sqrt{d}$-linéarisables. La démonstration repose sur le lemme \ref{lem4.1} ci-dessous. Commençons par 
introduire quelques définitions. 
On notera $S:=S_a+S_s$ la d\'ecomposition de Lebesgue d'un $(1,1)$-courant 
positif $S$. Celle-ci peut \^ etre d\'efinie \` a partir de la d\'ecomposition de Lebesgue des mesures
car un tel courant peut \^ etre consid\'er\'e comme une $(1,1)$-forme \` a coefficients mesures. Il est tr\`
es facile de voir que  
cette décomposition est 
unique et que les courants $S_a$, $S_s$ restent positifs. 
 Par contre la fermeture \'eventuelle de $S$ n'implique pas celle de $S_a$ ou de $S_s$. Nous noterons 
 $Supp(S)$ 
 le support de $S$ et $\sigma_S := S \wedge \omega^{k-1}$ sa mesure trace. On voit facilement que la 
 décomposition de Lebesgue de $\sigma_S$ est donnée par $\sigma_S=\sigma_{S_a}+\sigma_{S_s}$.  

\begin{lem}\label{lem4.1}
Soient $({\bf P}^k, f, \mu)$ un système $\sqrt{d}$-linéarisable, $S$ un courant positif de bidegré 
$(1,1)$ sur ${\bf P}^k$ tel que $f^{*}S=dS$ ($S$ n'est pas nécessairement fermé) et $\Omega$ un ouvert de 
${\bf P}^k$ chargé par $\mu$. 
\begin{itemize}
\item[1)]
 Si $S$ est absolument continu sur $\Omega$ ($S=S_a$) alors il existe une boule 
 $B(0,r)\subset {\bf C}^k$, un ouvert $\Omega' \subset \Omega$ chargé par $\mu$ et un biholomorphisme 
 $\Phi:B(0,r) \to \Omega' \subset \Omega$ tels que $\Phi^{*}S$ soit une forme différentielle à coefficients constants sur $B(0,r)$. 
\item[2)]  Supposons que $S$ dérive d'un potentiel \emph{psh} continu $v$ sur $\Omega$ $(S=dd^c v)$. Si 
$S_a$ est nul sur $\Omega$ alors $\mu(\Omega \cap Supp\;S) =0$.
\end{itemize}
\end{lem}

\underline{Démonstration du théorème \ref{TH}} :

Soit $\Omega$ un ouvert de ${\bf P}^k$ chargé par $\mu$. La première assertion du lemme \ref{lem4.1} 
appliquée à 
$T_a$ permet de supposer que dans de bonnes coordonnées, la restriction de $T_a$ à ${\Omega}$ est donnée 
par une forme
$H$ à coefficients constants. En particulier $T_a$ possède un
potentiel continu sur $\Omega$ et il en va donc de m\^ eme pour $T_s=T-T_a$ car $T$ est \`a potentiels
locaux continus. Ceci permet, sur 
$\Omega$, d'exprimer
$\mu$ sous la forme d'une somme de mesures positives obtenues comme
 produits extérieurs de $T_a$ et $T_s$ :
\begin{equation}\label{decmu}
\mu=T^k=\big(T_a+T_s\big)^k=T_a^k + \sum_{j=1}^k C_k^j\;\; T_s^j \wedge T_a^{k-j}.    
\end{equation}
Puisque $(T_s)_a$ est identiquement nul par définition, la seconde assertion du lemme \ref{lem4.1} montre 
que $\mu$ ne charge pas 
$\Omega\cap Supp\;T_s$ et
donc, au vu de (\ref{decmu}), la mesure $T_a^k$ n'est pas identiquement nulle sur $\Omega$. 
Autrement dit la forme $H$ n'est pas
dégénérée. Par ailleurs, puisque $\mu$ est absolument continue, chaque terme du second membre de 
(\ref{decmu}) doit, en tant que mesure positive, 
\^ etre absolument continue. En particulier,  la mesure singulière $T_s \wedge T_a^{k-1}$ est nulle sur 
$\Omega$. Or, $H$ étant strictement positive, celle-ci est équivalente à la mesure trace $\sigma_{T_s}$ de 
$T_s$.
 Le courant (positif) $T_s$ est donc nul sur $\Omega$ et $T$ coïncide sur cet ouvert avec une forme lisse 
 définie positive. L'endomorphisme $f$ est donc un exemple de Lattès comme cela est démontré dans \cite{BL}
 (voir le r\'esultat cit\'e dans l'introduction). \fin \\

\underline{Démonstration du lemme \ref{lem4.1}} :

Pour simplifier les notations, nous ne ferons pas figurer les cartes locales $\tau_x$ dans cette 
démonstration. Nous notons $\Lambda_n$ l'homothétie $(\sqrt{d})^{-n} \Id$.\\

1) Puisque $S$ est absolument continu sur $\Omega$, il est de la forme
$$ S=\frac{i}{2}\sum_{1\le p,q\le k} h_{p,q}(z)\;dz_p\wedge d\bar z_q \ \textrm{ où } 
\  h_{p,q} \in L^1(\Omega). $$

Soit $\cal M$ l'ensemble des points de $\Omega$ où toutes les fonctions $h_{p,q}$
sont continues en moyenne, c'est à dire : 
\[  \forall z \in {\cal M} \ , \ \lim_{r \to 0} {1 \over m(B(z,r))} \int_{B(z,r)} h_{p,q} (t) \, dm(t) = 
h_{p,q} (z).    \]
Puisque le syst\` eme est $\sqrt{d}$-lin\'earisable, $\mu$ est absolument continue par rapport \` a $m$
et l'ensemble $\cal M$ est de mesure totale pour $m$ et $\mu$. 
Notons $\cal R$ l'ensemble des points de $\Omega\cap Supp\;\mu$ où la suite $(f^n)_n$ est
linéarisable par des homothéties de rapport $(\sqrt{d})^{-n}$. Comme $\mu$ est absolument continue, la 
proposition
 \ref{pr4.1} nous assure que $\mu\big({\cal M}\cap{\cal R}\big)>0$. Soit alors 
$z \in {\cal M} \cap {\cal R}$ et posons $\Phi_n:=f^n \circ \Lambda_n$ (on identifie $z$ avec l'origine de  ${\bf C}^k$). Quitte à prendre une sous-suite, $\Phi_n(0)=f^n(z)$ reste dans $V \cap Supp\;\mu$ où $V$ est un voisinage de $z$ 
(cf Remarque \ref{remrec}) et la suite $(\Phi_n)_n$ converge vers un biholomorphisme 
$\Phi : B(0,\nu) \to \Omega' \subset \Omega$. Le support de $\mu$ étant fermé et invariant, on a 
$\Phi(0)\in \Omega \cap Supp\;\mu$ et donc  $\mu(\Omega')>0$. L'invariance de $S$ 
entra\^\i ne :    
$$\Phi_n^{*}S=\Lambda_n^{*} f^{n*}S=d^{n} \Lambda_n^{*} S=
\frac{i}{2}\sum_{1\le p,q\le k} h_{p,q}\circ \Lambda_n\; dz_p\wedge d\bar z_q.$$
Puisque $z \in {\cal M}$, on obtient 
$\Phi^{*}S=\frac{i}{2}\sum_{1\le p,q\le k} h_{p,q}(0)\;dz_p\wedge d\bar z_q$ par passage à la limite. \\

2) Supposons $\mu(\Omega \cap Supp S ) > 0$ et montrons que $S_a$ est non nul. Quitte à diminuer 
$\Omega$ on peut supposer que $S=dd^c v$ sur un voisinage $\tilde \Omega$ de $\overline{\Omega}$. 
Quitte à choisir une carte locale, $\tilde \Omega$ est un ouvert de ${\bf C}^k$. D'après la 
proposition \ref{pr4.1}, il existe ${\cal R} \subset \Omega \cap Supp\;S$ de $\mu$-mesure positive tel que 
pour tout point $z \in {\cal R}$, il existe une sous-suite  $\Phi_{n_j}:=f^{n_j} \circ (z + \Lambda_{n_j})$ convergeant uniformément sur $B(0,\nu(z))$ vers un biholomorphisme $\Phi$. On peut aussi supposer que $f^{n_j}(z) \in {\cal R}$ (cf remarque \ref{remrec}).
 
Observons tout d'abord qu'il suffit de  montrer que $\sigma_S$ possède une dérivée de Radon-Nykodym 
strictement positive en tout point 
$z$ de ${\cal R}$ : 
$$ \forall z  \in {\cal R} \ , \  \lim_n {1 \over d^{-kn}} \int_{B\big(z, \nu (\sqrt d)^{-n} \big)} 
S  \wedge \omega_0^{k-1} >0$$ 
où $\omega_0$ désigne la forme standard $\frac{i}{2} dd^c\vert\vert z\vert\vert^2$. 
En effet, comme $\mu({\cal R}) >0$, cette propriété montre que la mesure 
$\sigma_{S_a}$ (qui est égale à $(\sigma_S)_a$) n'est pas triviale sur $\Omega$ et il s'ensuit que 
le courant positif $S_a$ n'est pas nul.\\
V\'erifions \` a pr\'esent la stricte positivit\'e des d\'eriv\'ees. Notons que  
quitte à supprimer à ${\cal R}$ un ensemble de mesure de Lebesgue nulle 
(donc de $\mu$-mesure nulle), ces d\'eriv\'ees existent en tout point de $\cal R$.
Fixons donc $z \in {\cal R}$, et reprenons les applications $\Phi_{n_j}$ et $\Phi$ précédentes. Comme 
$\Phi(0) \in \overline{\Omega}$ on peut diminuer
$\nu$ de fa\c con à ce que les ouverts $\Phi_{n_j}\big(B(0,\nu)\big)$ et $\Phi\big(B(0,\nu)\big)$ soient 
contenus dans $\overline{\Omega}$. Puisque $f^{*}S=dS$, il vient :

$${1 \over d^{-k n_j}}  \int_{B \left( z, \nu  (\sqrt d)^{-n_j} \right) }  S  \wedge  \omega_0^{k-1}  =
{1 \over d^{-(k-1)n_j}}  \int_{z + \Lambda_{n_j} [B(\nu)]}   f^{{n_j}*}  S  \wedge  \omega_0^{k-1}=$$

$${1 \over d^{-(k-1)n_j}}  \int_{B(\nu)}  \Phi_{n_j}^{*}S\wedge\big(\Lambda_{n_j}^{*}\omega_0\big)^{k-1}=
\int_{B(\nu)}\Phi_{n_j}^{*}S\wedge\omega_0^{k-1}=
\int_{B(\nu)}dd^c(v\circ\Phi_{n_j})\wedge\omega_0^{k-1}$$
où $B(r)$ désigne la boule centrée en l'origine et de rayon $r$. Le théorème de convergence dominée 
entra\^\i ne alors :
$$\lim_j {1 \over d^{-k n_j}}  \int_{B(\nu \sqrt d)^{-n_j})}  S  \wedge   \omega_0^{k-1} \ge
\int_{B(\nu)}  dd^c(v\circ\Phi) \wedge   \omega_0^{k-1}=
\int_{B(\nu)} \Phi^{*}S         \wedge   \omega_0^{k-1}.$$
Cette dernière intégrale est bien strictement positive, car $\Phi(0) \in \overline{\cal R}\subset Supp\;S$. \fin

\section {Appendice}

Nous résumons ici la preuve du théorème du pluripotentiel utilisé dans la section \ref{diam}, ainsi que 
celle de l'estimation de la mesure présentée dans l'introduction.

\subsection{Un théorème de la théorie du pluripotentiel}

Il s'agit d'établir la version suivante d'un résultat d\^u à Briend-Duval \cite{BD} : \\

\noindent{\bf Théorème} :  
{\it Soit $S:=dd^c w$ un  $(1,1)$ courant positif fermé de potentiel $w$ continu sur le
polydisque $P(0,R)$ et 
$E\subset P(0,\frac{R}{2})$. On suppose que par tout $z \in E$ passe un disque holomorphe 
$\sigma_z : \Delta
\to {\bf C}^k$ de taille $l$ et qu'il existe une fonction $h_z$ harmonique sur $\Delta$
telle que 
$\vert w \circ \sigma_z - h_z \vert \le \epsilon$ sur $\Delta$. Alors il existe une constante $C(w)$ ne
dépendant que de $w$ telle que $S^k(E) \le C(w) \frac{k^2}{l^2}\epsilon$}.\\

Rappelons qu'un disque holomorphe $\sigma : \Delta \to {\bf C}^k$ passant par 
$z\in {\bf C}^k$ 
est dit de taille $l > 0$ si il est 
de la forme $\sigma(u)=z+l u. {v} +\beta(u)$ où ${v}$ est un vecteur unitaire de ${\bf C}^k$, 
$\beta(0)=0$ et $\vert\vert \beta \vert\vert\le \frac{l}{1000}$. \\

\underline{Démonstration} :
     
Soit $p_l$ la projection sur le $l$-ième axe de ${\bf C}^k$ et
$E_l:=\{z\in E / \vert\vert p_l(v_z) \vert\vert \ge \frac{1}{\sqrt{k}}\}$, de sorte que $E=\cup_{l=1,k}E_l$. 
Pour fixer les idées nous allons estimer $S^k(E_1)$.
A cet effet, on recouvre le polydisque $P(0,\frac{1}{2}R)$ par environ $N:=\frac{1}{4}
\frac{100 k}{l^2}$ ellipsoïdes contenus dans $P(0,R)$ de la forme ${\cal E}\big[
B(0,R)\big]$ où ${\cal E}(z_1,z')=\big(\frac{l}{10\sqrt{k}}z_1,z'\big)$.

Soit ${\cal E}$ l'un de ces ellipsoïdes. Puisque ${\cal E}$ est strictement pseudoconvexe, il existe une
fonction $\hat w$ $p.s.h$ maximale sur ${\cal E}$, continue sur $\overline {{\cal E}}$ et coïncidant avec
$w$ sur $b{\cal E}$.\\
Si $z\in {\cal E}\cap E_1$, on voit facilement que le disque $\sigma_z(\Delta)$ traverse ${\cal E}$, au 
sens
où la composante connexe  ${\cal C}$ de $\sigma_z^{-1}\big({\cal E}\cap\sigma_z(\Delta)\big)$ contenant 
l'origine
est relativement compacte dans $\Delta$. Un argument de principe du maximum montre que ${\cal C}$  
est simplement connexe. En l'exhaustant par des domaines à bord suffisamment régulier, on peut paramétrer
des disques holomorphes contenus dans $\cal E$ et dont le bord est arbitrairement proche de $b{\cal E}$.
Plus précisément, $\epsilon >0$ étant fixé, on trouve une
transformation conforme et continue jusqu'au bord $\psi : \overline{\Delta}\to 
\psi(\overline{\Delta})\subset{\cal
E}$ telle
que $\psi(0)=0$ et
$\vert \hat w - w \vert \le \epsilon$ sur $\sigma_z\circ\psi(b\Delta)$. Posons 
$\tilde{\sigma}_z:=\sigma_z \circ
\psi$ et notons $\tilde h$ la fonction harmonique sur $\Delta$ continue sur $\overline{\Delta}$ et
coïncidant avec $w\circ \tilde{\sigma}_z$ sur $b\Delta$. On a alors :
\begin{equation}\label{enca}
w(z)\le \hat w (z)\le \tilde h (0) +\epsilon
\end{equation}
la première inégalité provient de la maximalité de $\hat w$ sur $\cal E$ et la seconde du principe du
maximum appliqué à $\hat w\circ\tilde{\sigma}_z - \tilde  h$ (cette fonction coïncide avec 
$\hat w\circ\tilde{\sigma}_z - w\circ\tilde{\sigma}_z$ sur $b\Delta$).\\

Par hypothèse on a $h_z\circ\psi -\epsilon \le 
w\circ\sigma_z\circ \psi=w\circ\tilde{\sigma}_z
\le h_z\circ\psi +\epsilon$ sur $\overline \Delta$. 
On a donc aussi $h_z\circ\psi -\epsilon \le \tilde h \le h_z\circ\psi +\epsilon$ et il s'ensuit que :
\begin{equation}\label{maj2}
\vert w(z) - \tilde h(0)\vert \le 2\epsilon.
\end{equation}
Les inégalités (\ref{enca}) et (\ref{maj2}) montrent 
que :
$${\cal E}\cap E_1 \subset {\cal E}(w,\epsilon):=\{
z\in {\cal E} / 0\le \hat w(z) -w(z) \le 3\epsilon \}.$$ La majoration annoncée résulte alors 
immédiatement de
l'estimation suivante qui est au coeur de la démonstration de Briend-Duval et pour laquelle nous renvoyons 
à 
\cite{BD} ou \cite{S} page 180, Théorème A.10.2 :\\

{\it Il existe une constante} $C(w)>0$  {\it telle que} $(dd^c w)^k
\big[{\cal E}(w,\epsilon)\big]\le C(w) \epsilon$. \fin

\subsection{Estimation de la dimension}

Nous esquissons la preuve de l'estimation de la dimension en reprenant \emph{mutatis mutandis} les arguments développés par Binder et DeMarco \cite{BdM} dans le cas des endomorphismes polynomiaux de ${\bf C}^k$. \\

\noindent{\bf Théorème} :  
{\it Soit un système $({\bf P}^k, f, \mu)$ de degré $d$ et d'exposants 
$\lambda_1 \leq \cdots \leq \lambda_k$. La dimension de $\mu$ vérifie : 
$\dim(\mu) \leq 2(k-1) + {\log d \over \lambda_k}$ }.\\

Rappelons que la dimension est définie comme la borne inférieure des dimensions de Hausdorff des 
boréliens de mesure totale. Ce résultat montre que si la dimension de $\mu$ est égale à $2k$, alors 
tous les exposants de $\mu$ sont minimaux, égaux à $\log \sqrt{d}$.\\

\underline{Démonstration} :

\noindent Il s'agit d'exhiber pour tout $\epsilon > 0$ un borélien $Y$ de mesure totale vérifiant :
\begin{equation}\label{dim1}
\dim_H (Y) \leq 2(k-1)  + {\log d \over \lambda_k} +  {2k \over \lambda_k}\epsilon. 
\end{equation}
Soit  $\widehat A$ l'ensemble des points 
$\hat x = (x_{n})_{n \geq 0}$ de $\widehat \PP$ vérifiant pour tout $n \geq 0$ : 
\[  B(x_{-n} , {r_0 \over \kappa_0}e^{-n(\lambda_k + \epsilon)} )  \subset f^{-n}_{\hat x} [ B(x_0 ,r_0) ] 
\ \textrm{ et } \ m \left(  f^{-n}_{\hat x} [ B(x_0 ,r_0) ] \right) 
\leq \kappa_0  e^{-2n(\lambda_1+ \cdots + \lambda_k) + n \epsilon} .\]
On rappelle que $m$ désigne la mesure volume standard sur ${\bf P}^k$. On vérifie que si $\kappa_0$ est 
assez grand et $r_0$ assez petit, alors $\hat \mu (\widehat A)>0$ (cf \cite{BdM}, lemme 2).\\
 Soit $\widehat {A_n} := \hat f ^{-n} \widehat A$. 
La mesure $\hat \mu$ étant ergodique, le théorème de Birkhoff entra\^ \i ne que 
$\widehat Y := \limsup_n \widehat {A_n}$ est de mesure totale. 

On pose alors $Y := \pi_0 (\widehat Y)$ et $A_n := \pi_0 (\widehat {A_n})$, de sorte que 
$Y$ est aussi de mesure totale et est contenu dans $\limsup_n {A_n}$. Estimer la dimension de Hausdorff de 
$Y$ revient à estimer celle des ensembles $A_n$, pour $n$ assez grand. Par définition de $\widehat A$, 
tout point $y$ de $A_n$ vérifie :
\begin{enumerate}
\item
$f^n$ admet une branche inverse $g_n$ sur $B(f^n (y) ,r_0)$, telle que $g_n (f^n(y)) = y$

\item
La boule $B(y , {r_0 \over \kappa_0} e^{-n(\lambda_k + \epsilon)} )$ contient 
${\cal P} := g_n [ B(f^n (y) ,r_0) ]$

\item
$m ({\cal P} )  \leq k_0 e^{-2n(\lambda_1+ \cdots + \lambda_k) + n \epsilon}$.
 
\end{enumerate}

Il découle de ces propriétés que $A_n$ est recouvert par une famille 
$({\cal P}_i)_{i \in I}$ d'ouverts du type ${\cal P}$ dont le cardinal est de l'ordre de
$d^{kn}$. Pour le voir, il suffit de recouvrir $\overline{A_0}$ par un nombre fini de boules
$B(x_{i_0},\frac{1}{4}r_0)$ puis d'observer que tout $y\in A_n$ est dans 
$g_n\big[B(x_{i_0},\frac{1}{2}r_0)\big]$ dès lors que $f^n(y)\in B(x_{i_0},\frac{1}{4}r_0)$. 
D'après le point 3, 
le volume de la réunion des ${\cal P}_i$ n'excède pas 
$d^{kn} e^{-2n(\lambda_1+ \cdots + \lambda_k) + n \epsilon}$. 

Considérons à présent un recouvrement $({\cal M}_j)_{j \in J}$ de $A_n$ par des sous-ensembles de 
diamètre ${r_0 \over 100  \kappa_0} e^{-n(\lambda_k + \epsilon)}$ provenant d'un maillage de ${\bf P}^k$. 
D'après le point 2, 
un sous-ensemble ${\cal M}_j$ intersectant $A_n$ est nécessairement contenu dans 
$\cup_{i \in I} {\cal P}_i$. On a donc :   
\[ Card(J) \leq {  m ( \cup_{i \in I} {\cal P}_i )  \over m({\cal M}_j) } 
\lesssim { d^{kn} e^{-2n(\lambda_1+ \cdots + \lambda_k) + n \epsilon}   
\over \left(e^{-n(\lambda_k + \epsilon)} \right)^{2k} } .\]
En minorant les exposants $\lambda_1, \cdots, \lambda_{k-1}$ par $\log \sqrt{d}$, 
on obtient :
\[ Card(J) \lesssim  d^n  e^{n \left(  [2(k-1)\lambda_k] + (2k+1)\epsilon \right) }  .\]
Il s'ensuit que la mesure de Hausdorff de $A_n$, de dimension 
$l_\epsilon = 2(k-1) + \log d / \lambda_k + 2k \epsilon /\lambda_k$, est minorée par $e^{-n \epsilon}$
pour $n$ assez grand. La  $l_\epsilon$-mesure de Hausdorff de $Y \subset \limsup_n A_n$ est 
donc finie pour tout $\epsilon >0$. Cela termine la démonstration, car $Y$ est un borélien de mesure totale. \fin

\bigskip
\bigskip
\bigskip

{\footnotesize F. Berteloot}\\
{\footnotesize Universit\'e P. Sabatier, Toulouse III}\\
{\footnotesize Lab. Emile Picard, Bat. 1R2, UMR 5580}\\
{\footnotesize 118, route de Narbonne }\\
{\footnotesize 31062 Toulouse Cedex France}\\
{\footnotesize berteloo@picard.ups-tlse.fr}\\

\bigskip
\bigskip

{\footnotesize C. Dupont}\\
{\footnotesize Universit\'e Paris-Sud}\\
{\footnotesize Mathématique, Bat. 425, UMR 8628}\\
{\footnotesize 91405 Orsay, France}\\
{\footnotesize christophe.dupont@math.u-psud.fr}\\

\end{document}